\documentclass[11pt]{amsart}
\usepackage{amssymb,hyperref}
\usepackage{refcheck}

\newtheorem*{theorem}{Theorem}

\newcommand{\bcon}{\begin{conjecture}}
\newcommand{\econ}{\end{conjecture}}
\newcommand{\bcor}{\begin{corollary}}
\newcommand{\ecor}{\end{corollary}}
\newcommand{\bdf}{\begin{definition}}
\newcommand{\edf}{\end{definition}}
\newcommand{\beq}{\begin{equation}}
\newcommand{\eeq}{\end{equation}}
\newcommand{\bexa}{\begin{example}}
\newcommand{\eexa}{\end{example}}
\newcommand{\bexe}{\begin{exercise}}
\newcommand{\eexe}{\end{exercise}}
\newcommand{\bfac}{\begin{fact}}
\newcommand{\efac}{\end{fact}}
\newcommand{\bite}{\begin{itemize}}
\newcommand{\eite}{\end{itemize}}
\newcommand{\blem}{\begin{lemma}}
\newcommand{\elem}{\end{lemma}}
\newcommand{\bprb}{\begin{problem}}
\newcommand{\eprb}{\end{problem}}
\newcommand{\bpro}{\begin{proposition}}
\newcommand{\epro}{\end{proposition}}

\newcommand{\bque}{\begin{question}}
\newcommand{\eque}{\end{question}}
\newcommand{\brem}{\begin{remark}}
\newcommand{\erem}{\end{remark}}
\newcommand{\bthm}{\begin{theorem}}
\newcommand{\ethm}{\end{theorem}}
\newcommand{\bmat}{\begin{matrix}}
\newcommand{\emat}{\end{matrix}}

\newcommand{\bpr}{\begin{proof}}
\newcommand{\epr}{\end{proof}}

\newcommand{\comment}[1]{\,}

\title{Verification Of The Jones Unknot Conjecture Up To 23 Crossings}

\author{Robert E. Tuzun, Adam S. Sikora}

\begin{document}

\thispagestyle{empty}

\address{244 Math Bldg, University at Buffalo, SUNY, Buffalo, NY 14260}
\email{retuzun@buffalo.edu, asikora@buffalo.edu}

\pagestyle{myheadings}

\maketitle


The Jones conjecture states that the Jones polynomial distinguishes all non-trivial knots from the trivial one. Two years ago we confirmed that conjecture for all knots with diagrams up to 22 crossings. We described the details of that project in \cite{TS}. The purpose of this note is to announce an extension of this result:

\begin{theorem}
The Jones polynomial distinguishes all knots with diagrams up to 23 crossings and all algebraic knots with diagrams up to 24 crossings from the unknot. 
\end{theorem}

The method of this verification was that of \cite{TS} with further optimizations implemented. It involved the following steps:
\begin{enumerate} 
\item Generation of appropriate knot diagrams of 23 crossings, by (a) inserting ``algebraically-trivializable'' algebraic tangles into Conway polyhedra and by (b) considering closures of algebraic tangles.
\item Elimination of diagrams allowing a pass move which reduces either the number of crossings or the number of vertices in the corresponding Conway polyhedron.
\item Computation of the determinants of the remaining diagrams.
\item Computation of the Kauffman bracket polynomials of determinant $1$ diagrams by  the divide-and-conquer method.
\item  Computation of the knot group for the diagrams with monomial Kauffman brackets using the computer program SnapPy \cite{Sn}.
\item Checking the diagrams with ``non-trivial'' (non-cyclic) presentation of the knot group by other methods.
\end{enumerate}

The verification for 23 crossings required testing 10,675,582,685,989 non-algebraic knot diagrams and 43,356,077,900 algebraic knot diagrams, for a total of
10,718,938,763,889 knot diagrams -- approximately 9 times the number of knot
diagrams for 22 crossings.
For 24 crossings, 185,317,928,640 algebraic knot diagrams were tested.

Computations for 23 crossings were performed on 8-core Intel Xeon L5520 processors operated by the Center for Computational Research at the University at Buffalo for the total of 2744 core-days of wall-clock time (7.5 core-years).

%

\end{document}